\documentclass[10pt]{article}   
\usepackage{amssymb,amscd,latexsym}   
\usepackage{amsmath}
\usepackage{amsthm}
\newcommand{\rar}{\rightarrow}
\newcommand{\lar}{\longrightarrow}
\newcommand{\llar}{-\kern-5pt-\kern-5pt\longrightarrow}

\newcommand{\dashonto}{\dasharrow{\kern-15pt}\rightarrow}

\newtheorem{Theorem}{Theorem}[section]

\newtheorem{Proposition}[Theorem]{Proposition}
\newtheorem{Remark}[Theorem]{Remark}
\newtheorem{Example}[Theorem]{Example}
\newtheorem{Conjecture}[Theorem]{Conjecture}
\newtheorem{Definition}[Theorem]{Definition}

\def\sqr#1#2{{\vcenter{\hrule height.#2pt
        \hbox{\vrule width.#2pt height#1pt \kern#1pt
            \vrule width.#2pt}
        \hrule height.#2pt}}}
\def\phi{\varphi}
\def\demo{\noindent{\sc Proof. }}
\def\square{\mathchoice\sqr64\sqr64\sqr{4}3\sqr{3}3}
\def\qed{\hspace*{\fill} $\square$}

\DeclareMathOperator{\depth}{depth}

\def\xx{{\bf x}}
\def\yy{{\bf y}}
\def\zz{{\bf z}}

\def\ff{{\bf f}}

\def\ff{{\bf f}}
\def\gg{{\bf g}}

\def\fm{{\mathfrak m}}

\def\Ree#1{{\mathcal R}(#1)}

\def\depth{{\rm depth}\,}

\def\restr{{\kern-1pt\restriction\kern-1pt}}

\def\NN{\mathbb N}

\def\pp{{\mathbb P}}

\begin{document}
\begin{center}
{\Large{\bf\sc A theorem about Cremona maps and  symbolic Rees algebras}}
\footnotetext{2010 AMS {\it Mathematics Subject
Classification}: 13A30, 13C15, 13D02, 13H10, 14E05, 14E07.}

\vspace{0.3in}

{\large\sc Barbara Costa}
\quad\quad
{\large\sc Zaqueu Ramos}
\quad\quad
 {\large\sc Aron  Simis}\footnote{Supported by a CNPq grant and a PVNS Fellowship from CAPES.}

\end{center}


\bigskip

\begin{abstract}

This work is about the structure of the symbolic Rees algebra of the base ideal of a Cremona map. We give sufficient conditions under which this algebra has the ``expected form'' in some sense. The main theorem in this regard seemingly covers all previous results on the subject so far. The proof relies heavily on a criterion of birationality and the use of the so-called inversion factor of a Cremona map.
One adds a pretty long selection of examples of plane and space Cremona maps tested against the conditions of the theorem, with special emphasis on Cohen--Macaulay base ideals.

\end{abstract}

\section*{Introduction}
\label{intro}

The first central notion in this work is that of a Cremona map.
A Cremona map is a birational map of projective space $\pp^n_k={\rm Proj}(k[x_0,\ldots,x_n])$ ``onto'' itself. As a rational map, it has an associated ideal $I\subset R:=k[x_0,\ldots,x_n]$ which is generated by a linear system defining the map -- we call this ideal the {\em base ideal} of the map.

The main goal is to understand the structure of the symbolic Rees algebra of the base ideal of a Cremona map.
The main result is Theorem~\ref{Cremona-symbolic}, a result that we state in the case of a Cremona map whose base $I$ is such that the depth of the ring $R/I$ is positive.
In the plane case this restriction means that $R/I$ is a Cohen--Macaulay ring -- for short, we will say that a Cremona map is Cohen--Macaulay if $R/I$ is Cohen--Macaulay ring. For instance, any plane de Jonqui\`eres map is Cohen--Macaulay and so is any plane Cremona map of degree $\leq 4$.
In higher dimensional projective space the depth condition is often not very restrictive, allowing to include many important examples of Cremona maps which are not Cohen--Macaulay.

The structure of the symbolic Rees algebra is a delicate piece of algebra with roots in the early discoveries of  Zariski and Nagata.
It would be nearly insane to list the subsequent  papers written on the subject as soon as it became commutative algebra grounds -- however, see \cite[Introduction]{Zaron} and the references there.

Quite generally, given an ideal $I\subset R$, we are interested in the {\em symbolic filtration}, where
$I_{\ell}=I^{({\ell})}:=\iota^{-1}(S^{-1}I^{\ell})$, with $S\subset R$ standing for  the complement of the union of the associated prime ideals of $R/I$ on $R$ and
$\iota:R\rar S^{-1}R$ denotes the canonical homomorphism of fractions.
Clearly, this is exactly the intersection of the primary components of the ordinary power $I^{\ell}$ relative to the associated primes of $R/I$.
The associated Rees algebra, the {\em symbolic Rees algebra} $\mathcal{R}^{(I)}$, contains as graded $R$-subalgebra the Rees algebra
 $\mathcal{R}(I)$ associated to the filtration of the ordinary powers of $I$. Typically, the former will be described as an algebra extension of the latter.

We observe that for a radical ideal $I=P_1\cap\cdots\cap P_m\subset R$, one has simply
$I^{({\ell})}=P_1^{({\ell})}\cap\cdots\cap P_m^{({\ell})}$ and, more particularly, if each $P_j$ is a complete intersection
then $I^{({\ell})}=P_1^{\ell}\cap\cdots\cap P_m^{\ell}$.

There are more general definitions, of which the present notion is a particular case. In \cite[Section 3]{HHT} a relative concept is established where one is given another ideal $J$ and takes the saturation of a power $I^{\ell}$ by $J$.
In particular, this procedure allows to move along the primary components of the power of $I$ of increasing codimension  and introduce a symbolic power as far out as desired.
For example, one has the symbolic power which is
obtained by taking $J$ to be the intersection of all associated primes of the sequence $\{R/I^r\}_{r\geq 1}$ excluding the minimal primes of $R/I$.
Clearly this is the same as taking the multiplicative subset $S\subset R$ to be the complement of the union of the minimal primes of $R/I$ (instead of all associated primes of $R/I$ as we are taking).
Alas, this definition stumbles on a slight difficulty at the outset, for $\ell=1$, as the ideal $I$ itself may have embedded primes -- this is not uncommon, specially in the realm of base ideals of birational maps.
As a result, one would face a strange situation in which the generators of degree one of the symbolic algebra would come from the minimal primary component of $I$ instead of  $I$ itself. This means that the properties of the ideal $I$ would give room to a larger, perhaps not sufficiently subtle, ideal.

In order to circumvent this undesirable situation by avoiding restricting a priori the class of ideals under consideration, would force us to  define $I^{(1)}:=I$, making it a special rule for the symbolic power of order $1$.
The impact would  be that $\mathcal{R}(I)\subset \mathcal{R}^{(I)}$  would provide the minimal generators of the symbolic algebra in degree $1$.
However, this is a rather artificial solution.
Therefore, we will stick to the first notion above.
Thus, whenever we talk about $I^{({\ell})}$ we mean in that sense.
At the end, of course, this is a matter of choice, but regrettably the choice has been often made without context clearness.

\smallskip

In the first section we digress on the preliminaries of three subjects: birational maps, symbolic powers and inversion factors. 
In the first of these subjects, we explain a criterion of birationality in terms of Rees algebras that will be crucial.
The second brings the idea of considering the symbolic Rees algebra as a certain ideal transform, an idea that goes back to Rees (if not to Zariski itself). Here we state the main facts of this line of thought that will be applied.
The third subject emphasizes the role of the inversion factor of a Cremona map for obtaining generators of symbolic powers of certain order. 
These preliminaries appeared in some form elsewhere (see \cite{AHA}, \cite{Zaron}).

The second section contains the main theorem and its proof drawing upon the preliminaries of the previous section.
The statement is to the effect that the symbolic Rees algebra has the expected form under the mix of an asymptotic behavior condition and a purely ideal theoretic hypothesis.
The expected form is the extended Rees algebra of the base ideal of the inverse map to the given Cremona map.
Thus, there is no expected form of  the symbolic Rees algebra for the base ideal of an arbitrary rational map, it only making sense under the assumption of a Cremona map.

Section 3 is the longest. Here we give a selection of examples of Cremona maps with the purpose to applying the main theorem.
While this theorem implies that the symbolic Rees algebra has the expected form under those hypotheses, this form can still be proved valid under milder conditions. A purpose of these examples is to discuss these exceptions as well.
While entering this breach, we proved that the symbolic Rees algebra of a Cohen--Macaulay plane Cremona map has the expected form, although its base ideal typically fails for one of the conditions of the main theorem (Proposition~\ref{plane_C-M}).
The proof relies on a basic specialization result employed in \cite{ST}.

Tha last section is a brief experiment on a certain curious birational representation of $\pp^2$ in $\pp^3$.
One would think that such representations are geometrically totally understood, nevertheless their symbolic nature is far off reach.

\section{Preliminaries}

\subsection{Recap of birational maps}
\label{sec:1.1}

Our reference for the basics in this part is \cite{bir2003}, which contains enough of the introductory
material in the form we use here (see also \cite{AHA} for a more general overview).

Let $k$ denote an arbitrary infinite field which will be assumed to be algebraically closed for
the geometric purpose.
A rational map $\mathfrak{F}:\pp^n\dasharrow \pp^m$ is defined by $m+1$ forms $\ff=\{f_0,\ldots, f_m\}
\subset R:=k[\xx]=k[x_0,\ldots,x_n]$ of the same degree $d\geq 1$, not all null.
We often write $\mathfrak{F}=(f_0:\cdots :f_m)$ to underscore the projective setup.
Any rational map can without lost of generality be brought to satisfy  the condition
that $\gcd\{f_0,\cdots ,f_m\}=1$ (in the geometric terminology, $\mathfrak{F}$ {\em has no fixed part}).
Since our focus is on Cremona maps, we will  assume this condition throughout.
The common degree $d$ of the $f_j$ is the
{\em degree} of $\mathfrak{F}$ and
the ideal $I_{\mathfrak{F}}:=(f_0,\ldots,f_m)\subset R$ is called the {\em base ideal} of $\mathfrak{F}$.

The {\em image} of $\mathfrak{F}$ is the projective subvariety $W\subset \pp^m$ whose homogeneous
coordinate ring is the $k$-subalgebra $k[\ff]\subset R$ after degree renormalization.
Write $S:=k[\ff]\simeq k[\yy]/I(W)$, where $I(W)\subset k[\yy]=k[y_0,\ldots,y_m]$ is the homogeneous defining ideal
of the image in the embedding $W\subset \pp^m$.

We say that $\mathfrak{F}$ is {\em birational onto the image} if there is a rational map
backwards $\pp^m\dasharrow \pp^n$ such that the residue classes $\ff'=\{f'_0,\ldots, f'_n\}
\subset S$ of a set of defining coordinates do not simultaneously vanish and satisfy the
relations
\begin{equation}\label{birational_rule}\nonumber
(\ff'_0(\ff):\cdots :\ff'_n(\ff))=(x_0:\cdots :x_n), \;
(\ff_0(\ff'):\cdots :\ff_m(\ff'))\equiv (y_0:\cdots :y_m)\pmod {I(W)}
\end{equation}
There are of course other alternatives to introduce this notion in both algebraic or geometric language.
Since this work does not concern testing rational maps for birationality, we will not expand on those but as general references along the present line we suggest \cite{bir2003} and \cite{AHA}.

 A basic property (really, a characterization) of a birational map is that its graph and the graph of its inverse map can be suitably identified.
 We restate this from \cite[Proposition 2.1]{bir2003} so as to use it in subsequent sections.
 Quite generally, for an ideal $\mathfrak{a}\subset A$, we denote by $\mathcal{R}_A(\mathfrak{a})$ its Rees algebra $A[\mathfrak{a}t]\subset A[t]$.
 
 \begin{Proposition}\label{algcriterion} Let $X\subset \pp^n$ and $Y\subset \pp^m$ denote
 integral subvarieties of positive dimension and let $F\colon
 X\dasharrow \pp^m$ and $G\colon Y\dasharrow \pp^n$ stand for
 rational maps with  image $Y$ and $X$, respectively {\rm (}so that,
 in particular, $\dim X=\dim Y${\rm )}. Let $R=k[\xx]/\mathcal{I}(X)$ and
 $S=k[\yy]/\mathcal{I}(Y)$ denote the respective homogeneous coordinate rings
 of $X$ and $Y$. Fix sets of forms $\ff=\{f_0,\ldots,f_m\}\subset
 k[\xx]$ and $\gg=\{g_0,\ldots,g_n\}\subset k[\yy]$ whose respective
 residues in $R$ and $S$ are representatives of $F$ and $G$. The
 following are equivalent conditions:
 \begin{enumerate}
 \item[{\rm (i)}] $F$ and $G$ are inverse to each other.
 \item[{\rm (ii)}] The identity map of $k[\xx,\yy]/(\mathcal{I}(X),\mathcal{I}(Y))$ induces
 an isomorphism of bigraded $k$-algebras
 $${\cal R}_R\left(\frac{(\ff,\mathcal{I}(X))}{\mathcal{I}(X)}\right)
 \simeq {\cal R}_S\left(\frac{(\gg,\mathcal{I}(Y))}{\mathcal{I}(Y)}\right).$$
 \end{enumerate}
 \end{Proposition}

 The typical use we have in mind here is for the case of a Cremona map, in which case $R=k[\xx]$ and $S=k[\yy]$. 
 Set $I:=(\ff)\subset R$ and $I':=(\gg)\subset S$.
 In terms of the presentations $$\begin{array}{lcc}
 R[\yy]=k[\xx,\yy]=S[\xx] &\lar & S[I'u]\subset S[u]\\
\kern6pt \downarrow &&\\
 R[It]\subset R[t]&&
 \end{array}
 $$
 the above isomorphism $k[\xx,It]=R[It]\simeq S[I'u]=k[\yy,I'u]$ reads crosswise as $x_i\mapsto g_iu, y_i\mapsto f_it$.
Reading this isomorphism as identification takes place inside the full polynomial ring $k[\xx,\yy,t,u]$.
This innocent algebraic rule will show to be useful.

\medskip

Having information about the inverse map is often relevant.
For instance, for every representative $\ff'=\{f'_0,\ldots, f'_n\}$ of the inverse map, the structural congruence
(\ref{birational_rule}) involving the inverse map give a uniquely defined
form $D\in R$ such that $f'_i(f_0,\ldots,f_m)=x_iD$, for every $i=0,\ldots,n$.
We call $D$ the {\em source inversion factor} of $\mathfrak{F}$ associated to the given representative.
By  (\cite[Section 2]{AHA}), any such representative, as a vector, can be taken to be
 the transpose of a minimal generator of the syzygy module of the so-named weak Jacobian dual matrix.

\subsection{Recap of symbolic powers}
\label{sec:2}

In this part we recall  the Rees algebra of a (multiplicative) filtration
of ideals in a ring $R$.
Thus, let $\mathcal{F}=\{I_r\}_{r\in \NN}$ stand for such a filtration -- meaning that $I_r$ are ideals of $R$ such that
$I_r\subset I_{r-1}$ and $I_rI_s\subset I_{r+s}$, for any $r,s$.
The  Rees algebra of the filtration is the $\NN$-graded $R$-subalgebra
$$\mathcal{R}(\mathcal{F}):=R+I_1t+\cdots +I_rt^r+\cdots \subset R+Rt+\cdots +Rt^r+\cdots =R[t].$$
An important feature is that $\mathcal{R}(\mathcal{F})$ is a graded $R$-subalgebra of the standard graded polynomial ring
$R[t]$ with $R[t]_0=R$.

\begin{Definition}\rm We say that an element $f\in I_rt^r$ is an {\sc essential element}  {\sc of order} $r$ if it does not belong
to the graded piece of degree $r$ of the $R$-subalgebra $R[I_1t,\ldots,I_{r-1}t^{r-1}]\subset\mathcal{R}(\mathcal{F})$.
\end{Definition}
This amounts to say that

$$f\,\notin \sum_{
\begin{array}{c}
s_1+2s_2+\cdots +(r-1)s_{r-1}=r\\
0\leq s_j<r
\end{array}
}
R \,I_1^{s_1}I_2^{s_2}\cdots I_{r-1}^{s_{r-1}}= \sum_{1\leq s\leq r-1} R \, I_sI_{r-s}.$$

Thus, an essential element of order $r$ is a generator of $\mathcal{R}(\mathcal{F})$ in degree $r$
that does not come from elements in degree $r$ expressible as an $R$-combination of  elements of degree $\leq r-1$.

\smallskip

Given a fixed ideal $I\subset R$, we are interested in the {\em symbolic filtration}, where
$I_{\ell}=I^{({\ell})}:=\iota^{-1}(S^{-1}I^{\ell})$, with $S\subset R$ standing for  the complement of the union of the associated prime ideals of $R/I$ on $R$ and
$\iota:R\rar S^{-1}R$ denotes the canonical homomorphism of fractions.
Clearly, this is exactly the primary component of the ordinary power $I^{\ell}$ relative to the associated primes of $R/I$.

The associated Rees algebra is the {\em symbolic Rees algebra} of $I$.

To move deeper into the theory, one way is to use the nature of symbolic powers in terms of ideal transforms.
Namely, given an ideal $\mathfrak{a}\subset A$, where $A$ is a Noetherian domain with field of fractions $K$,
the {\em ideal transform} of $A$ relative to $\mathfrak{a}$ is the $A$-subalgebra $T_A(\mathfrak{a}):=A:_K\mathfrak{a}^{\infty}\subset K$.
We will draw on the following two fundamental principles:

\begin{Proposition}
{\rm (\cite[Proposition 7.1.4]{Wolmbook})}\label{fact1}
 If $C\subset T_A(\mathfrak{a})$ is a finitely generated $A$-subalgebra such that
${\rm depth}_{\mathfrak{a}C}(C)\geq 2$ then $C=T_A(\mathfrak{a})$.
\end{Proposition}

\begin{Proposition}
{\rm  (\cite[Proposition 7.2.6]{Wolmbook})}\label{fact2} 
Let $R$ be a Noetherian domain satisfying the condition $(S_2)$ of Serre and let $I\subset R$ be an ideal. Then
$$\mathcal{R}^{(I)}_R\simeq T_{\mathcal{R}(I)}(\mathfrak{a})\subset R[t]$$
as $R$-subalgebras of $R[t]$ for suitable choice of the ideal $\mathfrak{a}\subset R$.
\end{Proposition}

Here is an application of these principles:

\begin{Proposition}\label{idealtransform} Let $R=k[\xx]=k[x_0,\ldots, x_n]$  denote a standard graded polynomial ring over an infinite field $k$,
with irrelevant maximal ideal  $\fm:=(\xx)$.
Let $I\subset R$ stand for a homogeneous ideal with ${\rm depth}(R/I)>0$, satisfying the following properties:
\begin{enumerate}
\item[{\rm(i)}] For every $r\geq 1$, the $R$-module $I^{(\ell)}/I^{\ell}$ is either zero or $\fm$-primary.
\item[{\rm(ii)}]   If $C\subset \mathcal{R}^{(I)}_R$ is a finitely generated graded $R$-subalgebra 
containing the Rees algebra $\mathcal{R}_R(I)$ and such that ${\rm depth}_{\fm C}(C)\geq 2$, then $C=\mathcal{R}^{(I)}_R$.
 \end{enumerate}
\end{Proposition}
\demo
Since $R/I$ is assumed to have positive depth condition (i) implies that $\fm$ is the only other possible associated prime of any $R/I^{\ell}$  besides the minimal primes of $R/I$.
We may assume that some $I^{(\ell)}/I^{\ell}\neq 0$ otherwise there is nothing to prove.
Therefore $\fm$ is the only embedded associated prime of all powers,
i.e., we are actually taking the ordinary saturated Rees algebra.
This means that this algebra is the ideal transform on $C$ of the ideal $\fm$. Therefore, the result follows from Proposition~\ref{fact1}.
\qed

\subsection{The role of the inversion factor}
\label{sec:2.1}

In this part we explain how the two previous subsections come together.

Recall from Section~\ref{sec:1.1} the meaning of an inversion factor of a birational map
$\mathfrak{F}:\pp^n\dasharrow \pp^m$
onto its image.
For each representative $\gg$ of the inverse map $\mathfrak{F}^{-1}$, there is a uniquely defined source inversion factor.
By definition, these forms have the same degree $\deg(\mathfrak{F})\deg(\gg)-1$.

From now on, the base ring is a standard graded polynomial ring $R=k[\xx]=k[x_0,\ldots,x_n]$, with $\fm=(\xx)$.
The following property of the inversion factor plays a crucial role in the sequel.

\begin{Proposition}{\rm \cite[Proposition 1.4]{Zaron}}\label{inversionfactor_is_symbolic}
Let $\mathfrak{F}:\pp^n\dasharrow \pp^m$ be a rational map of degree $d\geq 2$ and let $I\subset R=k[\xx]$ denote its base ideal. Assume that $\mathfrak{F}$ is birational onto the image and let $D\subset R$ denote the source inversion factor relative to a given representative of the inverse map.

Suppose that ${\rm depth}(R/I)>0$. Then:
\begin{enumerate}
\item[{\rm (a)}] $D$ is an element of the symbolic power $I^{(d')}$,
where $d'$ is the degree of the coordinates of the representative.
In particular, $I^{(d')}\neq I^{d'}$.
\item[{\rm (b)}] If, moreover, $I^{(\ell)}=I^{\ell}$, $\ell\leq d'-1$, then $D$ is a essential symbolic element of order $d'$.
\item[{\rm (c)}] Moreover, if $I^{(d')}$ is generated in standard degree $\geq dd'-1$, where $d$ is the common degree of
 $\mathfrak{F}$, then $D$ is a homogeneous minimal generator of the symbolic Rees algebra.
\end{enumerate}
\end{Proposition}

The following elementary example illustrates the need for the triviality of the symbolic powers up to the required index.

\begin{Example}\rm
Let $I=(x_0x_1,x_1x_2,x_0x_2, x_2x_3,x_3x_4)$ -- this is the base ideal of a Cremona map of $\pp^4$.
By the method of \cite{CremonaMexico} (see also \cite{CostaSimis}),  the inverse map is defined by monomials of degree $3$
while the source inversion factor turns out to
be the monomial $x_0x_1x_2^2x_3$. Thus, the latter is an element in $I^{(3)}\setminus I^3$.
However, it is not a essential symbolic element of order $3$ because $x_0x_1x_2^2x_3=x_2x_3\cdot x_0x_1x_2\in II^{(2)}$.
\end{Example}
Note in this example the pattern: $x_0x_1x_2\in I^{(2)}$ is a homogeneous minimal generator of order $2$.

A large class of examples of the behavior in Proposition~\ref{inversionfactor_is_symbolic} is established in \cite{Zaqueu}
(see also \cite{Zaron}).
For convenience we give one instance of this class.

\begin{Example}\rm
Let  $I\subset R$ denote the ideal of maximal minors of a $4\times 3$ matrix whose entries are general
linear forms in $k[x_0,x_1,x_2,x_3]$.
Then $I$ is a radical ideal defining a Cremona map of $\pp^3$, satisfying $I^{(2)}=I^2$ and $I^{(3)}=(I^3, D)$, where $D$
is the source inversion factor. It can be shown that the degree of the inverse is again $3$, hence $\deg(D)=8$.
\end{Example}

\section{Symbolic algebras and Cremona maps}
\label{sec:2.2}

We are interested in the question whether in the case of the base ideal $I$ of a birational map its symbolic powers
are of a more special nature.

\subsection{Main theorem}

In this section we focus on the base ideal of a Cremona map to inquire about its symbolic powers.
We fix the notation: $R=k[X_1,\ldots,X_n]$ denotes a standard polynomial ring over an infinite field (often of characteristic zero).

The following theorem clarifies the impact of the Cremona context, unifying the  results of \cite[Corollary 2.4 (b)]{dual}, \cite[Theorem 2.14 (b)]{Zaron} and \cite[Proposition 2.5 (iii)]{Zaron2}.
We anticipate that, due to one of the crucial hypotheses,  the symbolic algebra turns out to be  the ordinary saturated Rees algebra (i.e., the symbolic algebra with respect to the irrelevant ideal).
However, the result is not a statement  about arbitrary such algebras.

\begin{Theorem}\label{Cremona-symbolic}
Let $\mathfrak{F}:\pp^n\dasharrow \pp^n$ stand for a Cremona map
with base ideal $I\subset R=k[\xx]=k[x_0,\ldots,x_n]$ such that ${\rm depth}(R/I)>0$.
Suppose that
\begin{enumerate}
\item[{\rm (i)}] For every $\ell\geq 1$ the module $I^{(\ell)}/I^{\ell}$ is either zero or $(\xx)$-primary
\item[{\rm (ii)}] The Rees algebra $\mathcal{R}(I)$ has the property $(S_2)$ of Serre.
\end{enumerate}
Then 
\begin{enumerate}
\item[{\rm (a)}] 
The symbolic algebra $\mathcal{R}^{(I)}=R[It,Dt^{d'}],$ where
$D$ denotes the  source inversion factor of $\mathfrak{F}$ and $d'$ is the degree of its inverse map.
\item[{\rm (b)}] If $\mathcal{I}\subset k[\xx,\yy]$ denotes the presentation ideal of $R[It]$ then the presentation ideal of $\mathcal{R}^{(I)}$ on $k[\xx,\yy,z]$ has the form
\begin{equation}\label{symbolic_ideal}
\left(\mathcal{I}, \{x_iz-g_i(\yy),0\leq i\leq n\}\right).
\end{equation}
\end{enumerate}
\end{Theorem}
\demo
(a) Using the notation of Proposition~\ref{algcriterion} and the comments after it, set $I=(f_0,\ldots,f_n)$ and $I'=(g_0,\ldots,g_n)\subset S$.
We will use the identification $R[It]=S[I'u]$ inside $k[\xx,\yy,t,u]$ as explained there.
Therefore, $x_i=g_iu, y_i=f_it$.
Now consider the element $Dt^{d'}\in R[t]=k[\xx,t]$, where $D$ is the source inversion factor of $\mathfrak{F}$ and $d'$ is the degree of the $g_i$'s.
By definition, say, $D=g_1(\ff)/x_1$.
Therefore 
\begin{eqnarray} 
Dt^{d'}&=&g_1(\ff)t^{d'}/x_1= g_1(\ff t)/x_1= g_1(\yy)/g_1(\yy)u\\
&=&u^{-1},
\end{eqnarray}
always by the above identification inside the field $k(\xx,\yy,t,u)$.
This gives $\mathcal{R}_R(I)[Dt^{d'}]= \mathcal{R}_S(I')[u^{-1}]$,
the latter being the extended Rees algebra of the ideal $I'\subset S$.

The rest of the argument mimicks the proof of  \cite[Corollary 2.4 (b)]{dual}, using the remaining hypotheses (i) and (ii) of the statement.
Namely, by (ii) $\mathcal{R}_S(I')=\mathcal{R}_R(I)$ satisfies $(S_2)$,
hence the associated graded ring ${\rm gr}_{I'}(S)$ satisfies $(S_1)$ by the usual exact sequences
$$\begin{array}{ccccccccc}
 0 &\rar & I'\mathcal{R}_S(I') &\lar & \mathcal{R}_S(I')& \lar & {\rm gr}_{I'}(S) &\rar & 0\\
&& \vert \wr  &&&&&&\\
0 &\rar & \mathcal{R}_S(I')_+ &\lar & \mathcal{R}_S(I')& \lar & S &\rar & 0
 \end{array}
 $$
 Since ${\rm gr}_{I'}(S)\simeq \mathcal{R}_S(I')[u^{-1}]/(u^{-1})$,
 with $u^{-1}$ a nonzerodivisor, it follows that $\mathcal{R}_S(I')[u^{-1}]$ satisfies $(S_2)$ and so does 
$\mathcal{R}_R(I)[Dt^{d'}]$ by the above identification.
This implies that the extended ideal $(\xx)\mathcal{R}_R(I)[Dt^{d'}]$ has depth at least $2$.
Finally use condition (i) to conclude by applying Proposition~\ref{idealtransform}.

(b) By the proof of (a),  $\mathcal{R}^{(I)}=\mathcal{R}_S(I')[u^{-1}]= \mathcal{R}_R(I)[u^{-1}]$.
Therefore, the assertion follows from the well-known presentation of the extended Rees algebra (see, e.g., \cite[Proposition 5.5.7]{HuSw}).
\qed

\begin{Remark}\label{relaxation}\rm
At the end of the proof above all we needed was that the extended ideal $(\xx)\mathcal{R}_R(I)[Dt^{d'}]$ have grade at least $2$.
Thus, condition (ii) can be relaxed in a particular situation. Of course, this weakening would be rather awkward to include in the statement of the theorem.
\end{Remark}

\section{Selected illustrative examples}

In this part we have chosen some examples whose discussion may extend beyond the usual. Besides conveying a good deal of working situations, the purpose  is also to hint at results concerning other typical hypotheses in the realm of Cremona transformations.

Whenever full application of  Theorem~\ref{Cremona-symbolic} is intended we will of course assume that $R/I$ has positive depth.
This is because we are taking the symbolic power relative to the ideal which is the intersection of all associated primes of the powers excluding the associated primes of $R/I$ itself.
Clearly, for this symbolic algebra, if $R/I$ has depth zero, it must reduce to the ordinary Rees algebra. 
Yet we often make a digression on the symbolic power obtained by taking instead the minimal primary parts of all powers.

\begin{Remark}\label{inverse_degree2}\rm
We emphasize that, regardless of any conditions -- such as (i) and (ii) in the statement of Theorem~\ref{Cremona-symbolic} -- for a Cremona map of $\pp^n$ whose inverse has degree $2$ the corresponding source inversion factor $D$ is always an essential generator of $I^{(2)}$ (this follows from Proposition~\ref{inversionfactor_is_symbolic} (b)).
In addition $\fm$ is an associated prime of all modules $I^{(\ell)}/I^{\ell}$ -- this is because, since $\fm$ drives $D$ inside $I^2$ then it drives $DI\subset I^{(2)}I\subset I^{(3)}$ into $I^3$, that is, $\fm=I^3:DI$, where $DI\not\subset I^3$.
Thus,  $\fm$ is an associated prime of $I^{(3)}/I^3$, and so on so forth.
\end{Remark}

\subsection{Cohen--Macaulay plane Cremona maps}

Since we are assuming that $R/I$ has positive depth, we are left in the plane case with Cremona maps whose base ideals are Cohen--Macaulay.
Moreover, a plane Cremona map and its inverse have equal degrees.
Finally, condition (i) will be automatic throughout.
This is because $\fm$ is an associated prime of all modules $I^{(\ell)}/I^{\ell}$. Since there cannot be any other associated primes besides $\fm$ and the minimal primes of $R/I$ then
$I^{(\ell)}=I^{\ell}: \fm^{\infty}$, for every $\ell\geq 1$, hence  $I^{(\ell)}/I^{\ell}$ either vanishes or is $\fm$-primary.

\begin{Example}\rm
The simplest case is of a Cremona map of degree $2$ as it is easily seen or well-known that the base ideal $I$ is Cohen--Macaulay and of linear type.
In particular, the Rees algebra of $I$ is a complete intersection -- hence condition (ii) of Theorem~\ref{Cremona-symbolic} is satisfied.
Thus, one can apply the main theorem.
\end{Example}

Perhaps more interesting:

\begin{Example}\rm 
Let $I\subset R=k[\xx]=k[x_0,x_1,x_2]$ denote the base ideal of a Cremona map of degree $\leq 4$.
\end{Example}
By \cite[Corollary 1.23]{HS}, any plane Cremona map of degree at most $4$ has a Cohen--Macaulay base ideal.
Moreover, any plane Cremona map of degree $\leq 3$ is a de Jonqui\`eres map (\cite[Proposition 2.9]{HS}) -- to be looked at in the next example.
Thus, we can focus on degree $4$.
There are only two types of homaloidal type in this degree as is known classically by a simple calculation of the equations of condition (see, e.g., \cite[The discussion after Proposition 2.9]{HS}; also \cite[Page 1, Footnote]{SeTy}).
Since one of these is a de Jonqui\`eres map, we are left with the only other possibility, namely, that $I$ is generated by the $2$-minors of a $3\times 2$ matrix whose columns have degree $2$.

Condition (ii) of the theorem is most likely satisfied in a strong form: the Rees algebra of $I$ is Cohen--Macaulay.
We give a proof under a certain hypothesis in the appendix at the end of the paper. We believe that this hypothesis is superfluous in the sense it can always be attained by a projective change of coordinates.

One might suspect that this might be the general behavior of plane Cremona maps that are not de Jonqui\`eres but whose base ideal is still saturated.
Unfortunately, already in degree $5$ this is false. 
In degree $5$ there is only one homaloidal type which is Cohen--Macaulay besides a de Jonqui\`eres map (see \cite[Theorem 2.14]{HS}).
It turns out that the Rees algebra of its base ideal does not satisfy $(S_2)$ -- perhaps somehow like in the de Jonqui\`eres case.
In degree $6$ the classification is within reach, but a lot more involved. Again here the Rees algebra fails to be Cohen--Macaulay.
The computation was done with {\em Macaulay}, but it gets hard pretty soon in higher degrees because often one has to take base points in some kind of general position.

\begin{Example}\rm
Let $I\subset R=k[x_0,x_1,x_2]$ denote the base ideal of a de Jonqui\`eres transformation. 
\end{Example}
One knows that $R/I$ is Cohen--Macaulay and by the above we may assume that the map is of degree $\geq 3$.
The only problem is condition (ii). 
In fact, the Rees algebra $\mathcal{R}(I)$ not only fails to be Cohen--Macaulay but it  does not satisfy $(S_2)$. To see this, we set $A:=R[y_0,y_1,y_2]=k[x_0,x_1,x_2,y_0,y_1,y_2]$.
It has been recently proved  (\cite[Corollary 1.6]{ST}) that the free graded resolution of $\mathcal{R}(I)$ has the shape
\begin{equation}\label{free_resolution_Jonq}
0\rar A^{d-3}\lar A^{2(d-2)}\lar A^{d} \lar A,
\end{equation}
where $d$ is the degree of the generators of $I$.
By a careful examination of the argument in the proof of this resolution and some other data collected in the details of the resolutions of the powers of $I$ as in \cite[Theorem 2.7]{HS}, one deduces that the entries of the tail matrix are either zero or the variables $x_0,x_1,y_1,y_2$ (for suitable reordering).
The prime ideal $\mathcal{Q}=(x_0,x_1,y_1,y_2)$ contains the defining ideal of 
$\mathcal{R}(I)$ on $A$ and localizing (\ref{free_resolution_Jonq}) at $\mathcal{Q}$ will give a minimal free resolution of same length.
Since we are now in a $4$-dimensional ring, the depth of $\mathcal{R}(I)_{\mathcal{Q}\mathcal{R}(I)}$ is $1$, while $\mathcal{Q}\mathcal{R}(I)$ is an ideal of codimension $2$.

\subsection{The symbolic algebra of a  Cohen--Macaulay plane Cremona map}

As we have seen, the main theorem is not applicable as such. However, note that the result of the theorem remains valid if one can prove directly that the extended ideal
\begin{equation}\label{extended_ideal}
\frac{\left(\fm, \mathcal{I}, \{x_iz-g_i(\yy),0\leq i\leq 2\}\right)}{\left(\mathcal{I}, \{x_iz-g_i(\yy),0\leq i\leq 2\}\right)}
\end{equation}
has grade $\geq 2$ in the ring $A[z]/\left(\mathcal{I}, \{x_iz-g_i(\yy),0\leq i\leq 2\}\right)$, where $\fm=(x_0,x_1,x_2)$.
As a side, note that, since the generators of $I$ are analytically independent. one has 
$$(\fm, \mathcal{I}, \{x_iz-g_i(\yy),0\leq i\leq 2\})=(\fm,\{g_i(\yy),0\leq i\leq 2\}),$$
an ideal of codimension $5$. Since the ideal $(\mathcal{I}, \{x_iz-g_i(\yy),0\leq i\leq 2\})$ has codimension $3$ (one more than the codimension of the defining ideal of the Rees algebra of $I$), we'd be done if codimension $\geq 2$ were the issue at stake.
Alas, we need to check the grade instead.

Luckily, nevertheless the needed argument goes through.
We state this in the following form:

\begin{Proposition}\label{plane_C-M}
Let $\mathfrak{F}:\pp^2\dasharrow \pp^2$ stand for a Cohen--Macaulay Cremona map. Then the conclusion of  Theorem~\ref{Cremona-symbolic} holds true.
\end{Proposition}
\demo
We proceed along the lines of \cite[Theorem 1.5 and Corollary 1.6]{ST}.
Namely, since $R/I$ has positive depth, up to a change of coordinates we may assume that $x_2$ is a regular element modulo $I$.

Letting $I\subset k[x_0,x_1,x_2]$ denote the base ideal of  $\mathfrak{F}$, set $J:=(I,x_2)/(x_2)$, which we identify with an ideal in $k[x_0,x_1]$ via the $k$-isomorphism 
$$k[x_0,x_1,x_2]/(x_2)\simeq k[x_0,x_1].$$
By \cite[Theorem 1.5 (i)]{ST}, the rational map $\eta:\pp^1\dasharrow \pp^2$ defined by the generators of $J$ is birational onto its image, a curve of degree $d$.

Let $\{g_0,g_1,g_2\}\in k[y_0,y_1,y_2]$ denote an ordered set of forms defining the inverse map to the given Cremona map.

\smallskip

{\sc Claim:} $g_2$ is up to a nonzero scalar the defining equation of the image curve of $\eta$.

\smallskip

To see this, we resort to the notation of \cite[Theorem 1.5]{ST}, whereby $\mathcal{J}=(\mathcal{A},E)$  where $\mathcal{A}\subset k[x_0,x_1,y_0,y_1,y_2]$ is such that $(\mathcal{I},x_2)=(\mathcal{A},x_2)$ and $E$ denotes the defining equation of the image curve and $\mathcal{J}$ denotes the Rees ideal of $J$.
On the other hand, as a consequence of the  Koszul-Hilbert Lemma of \cite[Proposition 2.1]{AHA}, one has  $x_0g_2-x_2g_1, x_1g_2-x_2g_3\in \mathcal{I}$. It follows that $x_0g_2,\;x_1g_2\in\mathcal{A}\subset \mathcal{I}$. Since $\mathcal{J}$ is a prime ideal, it readily obtains  $g_2\in\mathcal{J}.$  Then $g_2\in\mathcal{I}\cap k[y_0,y_1,y_2]=(E)$. Since  
 $\deg(g_2)=d=\deg(E),$ the claim follows. 
 
 \smallskip
 
 From the claim also follows that the inverse map to $\eta$ is defined by the residue classes of $g_0,g_1$ on the $k$-algebra $k[y_0,y_1, y_2]/(g_2)$.
 Moreover, the  extended Rees algebra of  this base ideal is  isomorphic to the algebra $T/x_2T$, where $$T:=\frac{k[x_0,x_1,x_2,y_0,y_1,y_2,z]}{(\mathcal{I},x_0z-g_0,x_1z-g_1,x_2z-g_2)},$$
 with $z$ an additional variable.
 Indeed, the extended Rees algebra is
 $$
 \frac{k[x_0,x_1,y_0,y_1,y_2,z]}{(\mathcal{J},x_0z-g_0,x_1z-g_1)}\simeq\frac{k[x_0,x_1,x_2,x_2,y_0,y_1,y_2,z]}{(\mathcal{J},x_0z-g_0,x_1z-g_1,x_2)}\simeq \frac{k[x_0,x_1,x_2,x_2,y_0,y_1,y_2,z]}{(\mathcal{I},x_0z-g_0,x_1z-g_1,x_2z-g_2,x_2)}\simeq T/x_2T,
 $$
 where we used the equalities 
 $$(\mathcal{J},x_2)=(\mathcal{A}, E,x_2)=(\mathcal{A}, g_2,x_2)=((\mathcal{A},x_2),g_2)=((\mathcal{I}
 ,x_2),g_2)=(\mathcal{I}, x_2z-g_2,x_2).$$
 
 \smallskip
 
 {\sc Claim:} The symbolic Rees algebra of $I$ is isomorphic to $T$.
 
 \smallskip
 
 As observed before, it suffices to show that the ideal $(x_0,x_1,x_2)T$ has grade at least $2$.
Since $T$ is a domain, $x_2$ is a regular element on $T$.
Now, we have seen that the defining ideal of $T/x_2T$ is also prime.
Since, e.g., $x_0$ does not belong to this ideal by a degree argument, then $\{x_2,x_0\}$ is a regular sequence on $T$, as was to be shown.
\qed

\subsection{Some Cremona maps in $\pp^3$}

\subsubsection{Polar Cremona maps}

\begin{Example}\rm 
Let $I\subset R=k[x_0,\ldots,x_3]$ denote the ideal of partials of the quartic $f=(x_1^2-x_0x_2)x_2x_3$, which is of codimension $2$.
Then $R/I$ is C-M.
\end{Example}

The polynomial $f$ is homaloidal -- i.e, $I$ defines a Cremona map -- since in this degree the Hilbert--Burch matrix of $I$ is linear and its Fitting ideals are easily seen to have enough codimension (to conclude use \cite[Example 2.4]{RuSi}). 
Moreover, the values of these codimensions also imply that  $I$ is of linear type, hence its Rees algebra is a complete intersection
(all this and more can be found in \cite{Trento}); in particular, condition (ii) of the theorem is fulfilled. 

Condition (i) of the same theorem fails here: a computation with {\em Macaulay} yields that the annihilator of $I^{(2)}/I^2$ is the codimension $3$ ideal $(x_1,x_2,x_0x_3)$.
 Therefore, the ideals $(x_0,x_1,x_2)$ and  $(x_1,x_2,x_3)$ are  associated primes of $I^{(2)}/I^2\subset R/I^2$, hence of $R/I^2$.
An additional calculation yields that these are the only associated primes of $R/I^2$ other than the minimal primes of $R/I$.
Moreover,  the set of associated primes of $R/I^3$ will be the stable value of ${\rm Ass}(R/I^{\ell})$ for $\ell >>0$.
This set is ${\rm Ass}(R/I)\cup \{(x_0,x_1,x_2),(x_1,x_2,x_3)\}\cup \{\fm\}$, where $\fm$ is the irrelevant maximal ideal.

Since $I$ is unmixed, to get a symbolic power $I^{\ell}$ we may either take its unmixed part or, alternatively, take the saturation $I^{\ell}:J^{\infty}$ with $J$ the intersection of all associated primes of all powers excluding the minimal primes of $R/I$. 
We choose the second method as it enhances the above calculation of associated primes of the powers.
Then $J=(x_1,x_2,x_0x_3)$.
To simplify the computation, we harmlessly replace this ideal by the single homogeneous element $x_1^2+x_2^2+x_0x_3\in J$ that does not belong to any minimal prime of $R/I$ (by an easy calculation, these are $(x_1,x_2), (x_2,x_3), (x_3, x_1^2-x_0x_2)$).

Writing $q:=x_1^2-x_0x_2$ and $c:=x_2^2x_3$, a computation with {\em Macaulay} then gives:

$$I^{(2)} \;\; \mbox{\rm has $2$ fresh generators, both of degree $5$}:\; cq,x_2x_3c$$
$$I^{(3)} \;\; \mbox{\rm has $1$ fresh generator, its degree is $8$}: \; x_1cq^2$$
$$I^{(4)} \;\; \mbox{\rm has $1$ fresh generator, its degree is $10$}: \; x_2cq^3$$

We note that neither $f\in I^{(2)}$ nor $f^2\in I^{(3)}$.
Moreover, the source inversion factor $D=c^2q\in II^{(2)}$, and hence is not a fresh generator of the symbolic algebra.
This should be compared to the behavior of the Cremona case of \cite{Zaron}, where the only fresh  generator is the source inversion factor in order $n=3$.

\medskip

{\sc Claim.}
The symbolic Rees algebra of $I$ is 
$$R[It, cq\,t^2,x_2x_3c\, t^2, x_1cq^2\, t^3, x_2cq^3\,t^4].$$

Letting $C$ denote this subalgebra of the full symbolic algebra, in order to apply Proposition~\ref{fact1} it  suffices to find an ideal $\mathfrak{a}\subset R$ whose ideal transform is the symbolic algebra. Such ideals can be obtained a priori if one knows the family of the associated primes of all powers of $I$ (cf. \cite[Proposition 7.2.6]{Wolmbook}.
In the present case the following ideal can be taken: $\mathfrak{a}=(x_1^2+x_2^2+x_0x_3, x_1x_2x_3)$, the first generator exactly as above and the second one a homogeneous element of $I$, together forming a regular sequence.

A harder computation with {\em Macaulay} will give the presentation ideal of $C$ on the polynomial ring $R[\yy,\zz]$ on two replicated sets of variables. The result is that $C$ is a Gorenstein algebra and that the ideal  $\mathfrak{a}C$ has codimension $2$ in $C$.
This proves the claim computationally, leaving little room for a more conceptual proof.

The interesting feature of this example is that if one takes instead the usual saturated powers (that is, with respect to the irrelevant maximal ideal) and form the corresponding symbolic algebra then it will have the same shape as in the main theorem, having the source inversion factor as essential generator of order $3$ (and the only essential generator in degree $\geq 2$).

\begin{Example}\rm
Let $I\subset R=k[x_0,\ldots,x_3]$ denote the ideal of partial derivatives of the determinant $f\in R$ of the $3\times 3$  degeneration of a Hankel matrix
\begin{equation}\label{subhankel}
\left(
  \begin{array}{ccc}
    x_0 & x_1 & x_2 \\
    x_1 & x_2 & x_3 \\
    x_2 & x_3 & 0 \\
  \end{array}
\right)
\end{equation}
as introduced in \cite[Section 4]{CRS}, where it is shown that $I$ defines a Cremona map.
\end{Example}
{\em Au contraire}, in this case $R/I$ has an embedded prime.
According to \cite[Theorem 5.6]{MaSi}, its associated primes are the (unique) minimal prime $P=(x_2,x_3)$ and the embedded prime $Q=(x_1,x_2,x_3)$.
However, one has $\depth R/I^2=1$, hence its associated primes are the same as of those of $R/I$. Therefore, we have $I^{(2)}=I^2$.

Moreover, we claim that for any $\ell\geq 1$, the only associated prime of $R/I^{\ell}$ of codimension $3$ is $Q$ again.

Indeed, one has  $I=(x_3^2,x_2x_3,-3x_2^2+2x_1x_3,x_1x_2-x_0x_3).$
Using the last and the first of these generators, one readily sees that $x_1^2x_2^2\in I$.
Assume that $Q'$ is a codimension $3$ prime ideal containing $\{x_2,x_3\}$, but not $x_1$.
Localizing we then get $x_2^2\in I_{Q'}$.
From the generators of $I$ this implies that $x_3\in I_{Q'}$, hence also $x_2\in I_{Q'}$.
Therefore,  $I_{Q'}=(x_2,x_3)_{Q'}=P_{Q'}$, generated by a regular sequence.
Clearly then, for any  $\ell\geq 1$, the only associated prime of
$R_{Q'}/(I_{Q'})^{\ell} $ is $P_{Q'}$.
Consequently, $Q'_{Q'}$ is not an associated prime of $R_{Q'}/(I_{Q'})^{\ell} $ and hence $Q'$ is not an associated prime of $R/I^{\ell}$.

On the other hand, the inverse of the Cremona map defined by $I$ is defined by forms of degree $3$, as it comes out of \cite[Theorem 5.8]{MaSi} (see \cite[Remark 4.6 (d)]{CRS} for an explanation).
(Actually, the forms of the inverse map are exactly the maximal minors of the $3\times 4$ Jacobian matrix with respect to $x_0,\ldots,x_3$ of the ``polynomialized'' linear syzygies of $I$.)
The corresponding source inversion factor $D$ is therefore of degree $5=2\cdot 3-1$. 
(Explicitly, up to a nonzero coefficient in $k$, it is $x_3^5$ -- this can readily be computed with {\em Macaulay}, but from \cite[Theorem 4.4 (iii) (d)]{CRS}) and \cite[Proposition 1.3]{Zaron} it follows at least that it is a multiple of $x_3^4$.)

Since previous symbolic powers are trivial, $D$ is an essential symbolic generator of order $3$. 
Furthermore, as we have seen, no module $R/I^{\ell}$ has associated primes other than those of $R/I$ and $\fm$.
It follows that $I^{(\ell)}=I^{\ell}:\fm^{\infty}$ and hence
 all quotients $I^{(\ell)}/I^{\ell}$ are primary to the irrelevant ideal. Thus, condition (i) of the main theorem is satisfied.

 Finally, an additional computation shows that the symmetric algebra of $I$ is Cohen--Macaulay; since $I$ is of linear type, condition (ii) of the theorem is satisfied.
 Therefore, this example is a clean application of that theorem.

\subsubsection{Cohen--Macaulay base ideal}

In this part we consider maps whose base ideal is Cohen--Macaulay.
Note that in this dimension this property is quite a bit stronger than being saturated.

\begin{Example}\rm
Let $I\subset k[x_0,x_1,x_2,x_3]$ denote the base ideal of a generalized de Jonqui\`eres map of $\pp^3$ of degree $\geq 3$ such that the underlying plane Cremona map has a (codimension $2$) Cohen--Macaulay base ideal.
\end{Example}
The algebraic background of this sort of map has been studied in \cite[Proposition 3.6 and Theorem 3.7]{NewCremona}, \cite[Section 3.1.1, Proposition 3.5]{PanSi}.
It has some remarkable features:

\smallskip

$\bullet$  $I$ is a codimension $2$ perfect ideal (i.e., $R/I$ is a Cohen--Macaulay ring of dimension $2$);

$\bullet$  $I=(qJ,f)$,
 where $J\subset k[x_0,x_1,x_2]$ is the base ideal  of the underlying plane Cremona map, and $q,f$ are $x_3$-monoids having no proper common factor, with $k[x_1,x_2]$-parts belonging to $J$;

$\bullet$ In degree $3$ the ideal $I$ is linearly presented and of linear type, while the corresponding map is an  involution.

\smallskip

The first two assertions and the ideal theoretic part of the third come out of the structure as described in the above references, while the claim about involution requires some additional work using that plane quadratic Cremona maps are themselves involutions (see \cite[Section 4]{SeTy} for an explicit description of this behavior).

We focus on the case of degree $3$.
Since  $I$ is of linear type, in particular its symmetric algebra is a complete intersection. Therefore, condition (ii) of the main theorem is satisfied.

Condition (i) is not satisfied, which can be seen in different ways.
Note that in this degree, $q$ is necessarily a linear form.
The minimal primes of $R/I$ are the minimal primes of $R/(q,f)$.
Since $q$ is linear, a typical minimal prime will be generated by $q$ and an irreducible factor of $f$ modulo $q$.
Therefore, embedded primes of the powers $R/I^{\ell}$ must contain some of these.
To compute a symbolic power $I^{(\ell)}$ one may saturate $I^{\ell}$ with respect to a form lying on each of these embedded primes, but not belonging to any of the minimal primes of $R/(q,f)$.

A computation with {\em Macaulay} shows that $R/I^2$ has embedded primes of codimension $3$ -- one for each minimal prime of $R/(q,f)$.
These will actually be all embedded primes of codimension $3$ of the powers.
The source inversion factor $D$ of the map is a square whose square root is an essential generator of $I^{(2)}$.
Then, in an indirect way,  the main theorem implies that condition (i) fails.

It is a fair guess that the symbolic algebra of $I$ is $R[It,\sqrt{D}t^2]$.

\begin{Example}\rm
Let $I\subset k[x_0,x_1,x_2,x_3]$ denote the base ideal of a Cremona map of $\pp^3$ of degree $3$ such that $R/I$ is Cohen--Macaulay of codimension $2$ and such that the inverse map is of degree $2$.
\end{Example}
In this situation one has necessarily:

\smallskip

$\bullet$ $I$ is linearly presented;

$\bullet$ $I$ is not of linear type;

$\bullet$ the base ideal of the inverse map is not Cohen--Macaulay.

\smallskip

The first is obvious since $I$ is generated by cubics.
The second assertion is a little subtler and follows from the first and the  hypothesis that the inverse map is quadratic (see \cite[Example 2.4]{RuSi}).
The third assertion is clear since four quadratic forms cannot be the maximal minors of a $4\times 3$ matrix.

\smallskip

Lat $\phi$ denote the presentation matrix of $I$ and $\fm$ the irrelevant maximal ideal.
In the present context, $I$ being of linear type means that  ${\rm cod}\,I_1(\phi)\geq 4$ and ${\rm cod}\,I_2(\phi)\geq 3$. (Since $I$ is necessarily linearly presented, the first bound means that $I_1(\phi)=\fm$.)
Thus, we are assuming that one of these bound is not attained.

We present two classical instances of this situation with slightly different behavior:

\begin{enumerate}
\item Noether Cremona
\item No-name Cremona
\end{enumerate}
 
The Noether Cremona is one of a series of Cremona maps devised by M. Noether (see \cite[Section 2.1]{RuSi}).
The base ideal is generated by the $3$-minors of the matrix
$$\left(
  \begin{array}{ccc}
    0 & -x_1 & -x_1 \\
    -x_0 & x_0 & x_1 \\
    x_0 & 0 & 0 \\
    x_2 & 0 & x_3
  \end{array}
\right)
$$
One easily sees that the $2$-minors are all contained in the codimension $2$ ideal $(x_0,x_1)$, hence $I$ is not of linear type.
Noether already had that the inverse map is of degree $2$, but an easy application of the criterion of \cite{AHA} gives explicit inverse coordinates.

By Remark~\ref{inverse_degree2}, the symbolic algebra has an essential generator $D$ of (degree $3$) and order $2$ and,
in addition, $\fm$ is an associated prime of all modules $I^{(\ell)}/I^{\ell}$.
Note that one cannot conclude as in the plane case that these modules are all $\fm$-primary for all $\ell\geq 1$ -- in principle, some could have associated primes of codimension $3$.
Actually, there is no a priori guarantee that $I^{(2)}$ itself does not have additional essential generators of degree $5$ or less.
However, a computation with {\em Macaulay} shows that this cannot happen, so condition (i) is satisfied.

The no-name Cremona has a base ideal the $3$-minors of the following matrix:
$$\left(
  \begin{array}{ccc}
    2x_0 & 0 & 0 \\
    x_1 & 2x_0 & 0 \\
    0 & 3/2x_1 & 2x_0 \\
    -x_3 & x_2 & x_1
  \end{array}
\right)
$$
Its inverse map is given by the polar map of the determinant of the matrix (\ref{subhankel}) studied earlier.
As the Noether Cremona, this one follows the same pattern as regards the symbolic algebra, satisfying the hypotheses of the main theorem.

\begin{Remark}\rm
Although symbolically very similar, not to mention the matrix appearances, the Noether and the no-name maps have algebraic-geometric distinct features. Thus, for example, Noether's has a unique point in its base set, while the no-name one has $4$ such points.
Besides, the inverses, both of degree $2$, are very different in that $R/I$ is equidimensional in the no-name case but not in Noether's.
\end{Remark}

\section{Elements of a birational representation of $\pp^2$ in $\pp^3$}
\label{2:4}

Here we seek to understand the behavior of the symbolic powers of the base ideal of a certain birational representation $\pp^2\dashonto W\subset \pp^3$.
Due to the complexity of the subject, we will restrict ourselves to the case where the ideal is a codimension $2$ Cohen--Macaulay ideal.

This section was inspired  by the curious example in \cite[Example 4.6]{SUV95}, whereas its objective is a modest tentative to extend the setup  in \cite[Section 2.4]{Zaron} to the case where the structural matrix is allowed to have one single column of arbitrary degree.

\subsection{The ideal}

For convenience we will change notation slightly, writing $R:=k[x_1,\ldots,x_n]$.
The basic template of this section is an $(n+1)\times n$ matrix
\begin{equation}\label{template_matrix}
\phi_{n,r}=\left(
  \begin{array}{cccc}
   \lambda_{1,1} & \cdots &  \lambda_{1,n-1} & \rho_1 \\
    \vdots &  \vdots & \vdots  & \vdots\\
     \lambda_{n+1,1} & \cdots &  \lambda_{n+1,n-1} & \rho_{n+1}
  \end{array}
\right)
\end{equation}
wher $\lambda_{i,j}$ is a linear form and $\rho_i$ is a form of degree $r\geq 1$.
Let $I\subset R$ denote the ideal generated by the $n$-minors of this matrix. We assume that it has codimension $\geq 2$, hence $R/I$ is a codimension $2$ Cohen--Macaulay ring generated by forms of degree $r+n-1$.

The homological behavior of $R/I$ has quite an involved numerology.
From its free resolution
$$0\rar R(-(r+n))^{n-1}\oplus R(-(2r+n-1))\lar R(-(r+n-1))^{n+1}\lar R$$
it is straightforward to compute its multiplicity $e_0(R/I)=r^2+(n-1)r+{n \choose 2}$.
From this, if the rational map defined by $I$ is birational onto the image we get a weak upper bound for the degree of the implicit equation ${\rm edeg}(I)$:
$$ {\rm edeg}(I)\leq (n+r-1)^{n-3}\left((r+n-1)^2-e_0(R/I)\right)= \frac{(n+r-1)^{n-3}(n-1)(2r+n-2)}{2}.$$
If $n=3$ this becomes an equality (see, e.g., \cite[Theorem 6.6 (a)]{ram2})).
Therefore, for $n=3$ one has ${\rm edeg}(I)=2r+1$, and this is the only case we will be dealing with.

\subsection{The Rees algebra}

The punch line as regards the linearly presented case is that, for $r\geq 2$, the Rees algebra $\mathcal{R}(I)$ is neither of the expected form nor Cohen--Macaulay.
We nevertheless have the following result that gives an approximation to the structure of its defining ideal.
We set $k[\xx]=k[x_1,x_2,x_3]$ and $k[\yy]=k[y_1,y_2,y_3,y_4]$.

\begin{Proposition}
Consider the matrix $\phi:=\phi_{3,r}$ as in {\rm (\ref{template_matrix})}, and let $I:=I_3(\phi)\subset k[\xx]$.
Assume that the $2$-minors of the linear part of $\phi$ is $(\xx)$-primary.
Then:
\begin{enumerate}
\item[{\rm (a)}] $I$ defines a birational map onto the image
and admits an inverse representative of degree $2\,${\rm ;}
\item[{\rm (b)}] For every $i=1,2,\ldots, r$, let $f_i$ denote the Sylvester form of $\{l_1,l_2,f_{i-1}\}$ with respect to $\{\xx\}$, where $l_1,l_2, f_0$ are the forms generating $I_1((\yy)\cdot \phi)$.
Then the defining ideal of $\mathcal{R}(I)$ is a minimal prime of the ideal
\begin{equation}
\left(I_1\left((\yy)\cdot \phi\right), f_1,f_2\ldots, f_{r}
\right),
\end{equation}
where $f_i$ has bidegree $(r-i,2i+1)$, for $i=1,2,\ldots, r\,$.
\end{enumerate}
\end{Proposition}
\demo
(a) Consider the forms $l_1,l_2\in k[\xx,\yy]$ obtained from the first two columns of $\phi$. Because of the standing assumption on the $2$-minors of these two columns the Jacobian  matrix $\Theta$ of $l_1,l_2$ with respect to $\xx$ has rank $2$ over $k[\yy]$, hence also modulo the implicit equation, since the latter has degree $2r+1\geq 3$ by a previous remark.
By the main criterion of \cite{AHA}, the map defined by $I$ is birational onto the image with inverse being defined by the $2$-minors of $\Theta$. 

\smallskip

(b) Let $\mathcal{J}\subset k[\xx,\yy]$ denote a presentation ideal of 
$\mathcal{R}(I)$. Clearly, $I_1((\yy)\cdot \phi))\subset \mathcal{J}$, since the former is a presentation ideal of the symmetric algebra of $I$.
By induction on $i$ it is straightforward to write the degrees of the iterated Sylvester forms as stated; clearly, we know that, by construction, these forms belong to $\mathcal{J}$.

Since the ideal thus far generated has codimension $\geq 3$,
we are through.

\begin{Conjecture}\rm
The ideal in item (b) of the above proposition is indeed the Rees ideal of $I$.
\end{Conjecture}
The shortest way to prove the conjecture would be to show that the ideal is prime.
But this seems to be the hardest as well.

\subsection{The symbolic algebra}

The pattern of the symbolic powers seems to be quite irregular.
First, for any $r\geq 1$, $I^{(2)}/I^2$ has at least one minimal generator given by the source inversion factor coming from an inverse representative of degree $2$, which has degree $2(r+2)-1=2r+3$.
But this seems to be the only common feature for all $r$.
Indeed, one has this initial irregular behavior:

$\bullet$ For $r=1$, $I^{(2)}/I^2$ is minimally generated by the three source inversion factors coming, respectively, from the three inverse representatives of degree $2$

$\bullet$ For $r=2$, $I^{(2)}/I^2$ is minimally generated by the source inversion factor coming from the (unique) inverse representative of degree $2$ and two more forms of the same degree as the source inversion factor

$\bullet$ For $r=3$, $I^{(2)}/I^2$ is minimally generated by the source inversion factor coming from the (unique) inverse representative of degree $2$ (hence is cyclic in this case)

$\bullet$ For $r=4$, $I^{(2)}/I^2$ is minimally generated by the source inversion factor coming from the (unique) inverse representative of degree $2$ and two more forms of degree $13$ (hence, larger than the initial degree of $I^2$)

\section{Appendix}

\begin{Proposition}
Let $I\subset R=k[x_0,x_1,x_2]$ denote the base ideal of a plane Cremona map of degree $4$.
AssumeThat:
\begin{enumerate}
\item[{\rm (i)}] The map is not a de Jonqui\`eres map
\item[{\rm (ii)}] No entry of the $3\times 2$ syzygy matrix of $I$ 
admits a nonzero pure power term.
\end{enumerate}
Then the Rees algebra $\Ree I$ is a normal Cohen--Macaulay ring.
\end{Proposition}
\demo
Note that the inverse map is not a de Jonqui\`eres map either and is of degree $4$, hence its base ideal $I'\subset S=k[\yy]= k[y_0,y_1,y_2]$ is Cohen--Macaulay, and the syzygy matrix of $I'$ is also formed with columns of degree $2$.
On the other hand, one has an equality of the presentation ideals of the respective Rees algebras of $I$ and $I'$ on the ring $k[\xx,\yy]$ -- call it $\mathcal J$. This implies that among the generators of $\mathcal J$ one has at the outset four forms: two of bidegree $(2,1)$ and two of bidegree $(1,2)$.
The result will be that the ideal generated by these forms is prime -- in which case it must coincide with the presentation ideal ${\mathcal I}\subset k[\xx,\yy]$ -- and the forms are the maximal minors of a $4\times 3$ matrix.

In order to prove it we proceed as follows. 
Without loss of generality, by assumption (ii) the syzygy matrix of $I$ has the form 

$$\varphi=\left(\begin{array}{cc}
a_{11}x_0x_1+b_{11}x_0x_2+c_{11}x_1x_2&a_{12}x_0x_1+b_{12}x_0x_2+c_{12}x_1x_2\\
a_{21}x_0x_1+b_{21}x_0x_2+c_{21}x_1x_2&a_{22}x_0x_1+b_{22}x_0x_2+c_{22}x_1x_2\\
a_{31}x_0x_1+b_{31}x_0x_2+c_{31}x_1x_2&a_{32}x_0x_1+b_{32}x_0x_2+c_{32}x_1x_2
\end{array}
\right)$$

The  presentation of the symmetric algebra of $I$ can therefore be generated by the following two forms of bidegree $(2,1)$:
 $$\left(\sum a_{i1}y_i\right)x_0x_1+\left(\sum b_{i1}y_i\right)x_0x_2+\left(\sum c_{i1}y_i\right)x_1x_2=\det\left(\begin{array}{ccc}\sum a_{i1}y_i&\sum b_{i1}y_i&\sum c_{i1}y_i\\
 -x_2&x_1&0\\
 -x_2&0&x_0\end{array}\right)=:\Delta_2$$
 
 and

  $$\left(\sum a_{i2}y_i\right)x_0x_1+\left(\sum b_{i2}y_i\right)x_0x_2+\left(\sum c_{i2}y_i\right)x_1x_2=\det\left(\begin{array}{ccc}\sum a_{i2}y_i&\sum b_{i2}y_i&\sum c_{i2}y_i\\
 -x_2&x_1&0\\
 -x_2&0&x_0\end{array}\right)=:\Delta_1.$$
 
 Consider the matrix
 $$B=\left(\begin{array}{ccc}
 \sum a_{i1}y_i&\sum b_{i1}y_i&\sum c_{i1}y_i\\
 \sum a_{i2}y_i&\sum b_{i2}y_i&\sum c_{i2}y_i\\
  -x_2&x_1&0\\
  -x_2&0&x_0\end{array}\right),$$
  of which the above two determinants are maximal minors.
  
  We claim the remaining two maximal minors also belong to $\mathcal{J}$.
  Indeed, Laplace expansion along, respectively, the third and fourth rows of $B$ yields the minors

\begin{equation}\label{Delta_3}
\Delta_3:=q_3(\yy)x_0-q_1(\yy)x_2
\end{equation}

and 

\begin{equation}\label{Delta_4}
\Delta_4:=-q_2(\yy)x_1-q_1(\yy)x_2
\end{equation}
where
$$q_1(\yy)=\left(\sum b_{i1}y_i\right)\left(\sum c_{i2}y_i\right)-\left(\sum b_{i2}y_i\right)\left(\sum c_{i1}y_i\right)$$

$$q_2(\yy) =\left(\sum a_{i1}y_i\right)\left(\sum c_{i2}y_i\right)-\left(\sum a_{i2}y_i\right)\left(\sum c_{i1}y_i\right)$$ 

$$q_3(\yy)=\left(\sum a_{i1}y_i\right)\left(\sum b_{i2}y_i\right)-\left(\sum a_{i2}y_i\right)\left(\sum b_{i1}y_i\right)$$
are the $2$-minors of the first two rows of $B$.
On the other hand, applying Cramer along the second and third columns of $B$, respectively, we find the relations

\begin{equation}\label{eq3}
\left(\sum b_{i2}y_i\right)\Delta_2-\left(\sum b_{i1}y_i\right)\Delta_1=(q_3(\yy)x_0-q_1(\yy)x_2)x_1
\end{equation}
and

\begin{equation}\label{eq4}
\left(\sum c_{i2}y_i\right)\Delta_2-\left(\sum c_{i1}y_i\right)\Delta_1=(-q_2(\yy)x_1-q_1(\yy)x_2)x_0,
\end{equation}
This shows that the right hand sides are elements of $\mathcal{J}$.
But since the latter is a prime ideal not containing any $x$-variable, we deduce that $\Delta_3,\Delta_4\in \mathcal{J}$ as well.  
Since they are both forms of bidegree $(1,2)$, they must be the generators of the ideal of the symmetric algebra of the ideal $I'$
provided none of them vanishes.

\smallskip

We proceed to show that $I_3(B)=(\Delta_1,\Delta_2,\Delta_3,\Delta_4)\subset \mathcal{J}$ is an equality.
For it, we show namely:

\smallskip

$\bullet$ $I_3(B)$ has codimension $2$, hence $k[\xx,\yy]/I_3(B)$ is Cohen--Macaulay;

$\bullet$ $q_i\neq 0$, for $i=1,2,3$;

$\bullet$ $k[\xx,\yy]/I_3(B)$ is locally regular in codimension $1$.

\smallskip

From these follows that $k[\xx,\yy]/I_3(B)$ is normal, hence integral as $I_3(B)$ is a homogeneous ideal.
Since $I_3(B)$ is prime then the desired equality follows.

The first assertion is the fact that the symmetric algebra of $I$ has dimension $4$ because the heights of the Fitting ideals of $\phi$ are the expected ones. Therefore, $(\Delta_1,\Delta_2)$ has codimension $2$.

The second assertion is taken care by the observation that if, say, $q_1=0$ then $\Delta_3:=q_3(\yy)x_0, 
\Delta_4:=-q_2(\yy)x_1$, and hence $q_2$ and $q_3$ belong to $\mathcal{J}$, which, since they are forms in $\yy$, is only possible if they vanish.
Therefore, it suffices to show that one of them, say, $q_1$, is nonzero.
Suppose $q_1=0$, i.e.,  
\begin{equation}\label{eq1}
\left(\sum b_{i1}y_i\right)\left(\sum c_{i2}y_i\right)=\left(\sum b_{i2}y_i\right)\left(\sum c_{i1}y_i\right),
\end{equation}
Note that the factors on each hand side are nonzero linear forms  as otherwise the ideal $I$ would be contained in the ideal generated by one of the $\xx$-variables.
By a similar token, $\sum b_{i1}y_i$ does not divide $\sum b_{i2}y_i$ as otherwise an appropriate elementary column operation would lead to the same situation in which the ideal $I$ would be contained in the ideal generated by one of the $\xx$-variables.
Thus, we must conclude that 
$\sum b_{i1}y_i$ is a multiple of  $\sum c_{i1}y_i$.
Applying a similar reasoning, $\sum a_{i1}y_i$ is a multiple of  $\sum c_{i1}y_i$.
Substituting in $\phi$ and applying elementary column operations we obtaing a column with two zero entries. This would say that $I$ is generated by $2$ elements, an absurd.

\medskip

Finally, for the third above assertion, we check the Jacobian ideal of $I_3(B)$.
The transposed Jacobian matrix of $\{\Delta_1,\Delta_2,\Delta_3,\Delta_4\}$ has the form

 $$
\Theta =\left(\begin{array}{cc}
* & \phi'\\
\phi & *
\end{array}
\right)
$$
where 
$$
\phi'=\left(\begin{array}{cc}
q_3(\yy) & 0\\
0 & -q_2(\yy)\\
-q_1(\yy) & -q_1(\yy)
\end{array}
\right)$$
Since $\phi$ and $\phi'$ are in separate sets of variables, $I_2(\Theta)$ has height at least $4$ provided $I_2(\phi')$ has height $2$. 
But by the second assertion none of the entries of $\phi'$ vanishes.
This implies that $\phi'$ has rank $2$.
On the other hand, it is submatrix of the full weak Jacobian dual of $I$. By \cite[Theorem 2.18 (ii)]{AHA}, the inverse map is defined by the $2\times 2$ minors of $\phi'$. But since the inverse map is of degree $4$, the $\gcd$ of these $3$ minors is $1$. This means that $I_2(\phi')$ has height $2$.
\qed

\begin{Example}\rm
The following example appears in \cite[Example 2.1.8]{alberich}:
\end{Example}
$$I=(\,x_0x_1(x_1-x_2)(x_0-x_2),\, x_0x_2(x_1-x_2)(x_0-2x_1),\, x_1x_2(x_0-x_2)(x_0-2x_1)\,).
$$
As it turns out, $I$ coincides with the ideal generated by the degree $4$ part of the fat ideal of the six points with the following coordinates
$$\left(\begin{array}{cccccc}
1 & 0 & 0  & 1 & 2 & 2\\
0 & 1 & 0  & 1 & 1 & 1\\
0 & 0 & 1 & 1 & 1 & 2
\end{array}\right)
$$
with respective multiplicities $2,2,2,1,1,1,$ (actually, the fat ideal itself is generated in degree $4$).
The syzygy matrix of $I$ can be written as

$$\left(\begin{array}{cc}
-x_2(x_0-2x_1) & 0\\
x_1(x_0-x_2) & x_1(x_0-x_2)\\
0 & x_0(x_1-x_2)
\end{array}\right),
$$
therefore it satisfies condition (ii) of the proposition.

\begin{Remark}\rm
If one is willing to accept points with random coordinates then a computation with {\em Macaulay} gives that the Rees algebra of $\mathcal{J}$ is Cohen--Macaulay with the presentation ideal of the nature obtained in the proof of the proposition.
Moreover, replacing $3$ of the $6$ random points by the coordinate points and imposing that these have virtual multiplicity $2$ each reinstates condition (ii) of the proposition.
It would be nice giving a rigorous proof that, in general, if we apply a projective change of coordinates to have $3$ points with virtual multiplicity $2$ each to be the coordinate points, then condition (ii) holds. This would show that the Cohen--Macaulayness of the Rees algebra always holds.
\end{Remark}

\noindent {\bf Authors' addresses:}

\medskip

\noindent {\sc Barbara Costa}, Departamento de Matem\'atica\\ Universidade Federal Rural de Pernambuco \\
              52171-900, Recife, PE, Brazil\\
 {\em e-mail}:  binhamat@gmail.com \\

\noindent {\sc Zaqueu Ramos}, Departamento de Matem\'atica, CCET\\ Universidade Federal de Sergipe\\
49100-000 S\~ao Cristov\~ao, Sergipe, Brazil\\
{\em e-mail}: zaqueu.ramos@gmail.com\\

\noindent {\sc Aron Simis},  Departamento de Matem\'atica, CCEN\\ Universidade Federal
de Pernambuco\\
 50740-560 Recife, PE, Brazil.\\
{\em e-mail}:  aron@dmat.ufpe.br

\end{document}